\documentclass[11pt]{article}

\usepackage{amsmath,amssymb,amsfonts}
\usepackage{fullpage}
\usepackage[all]{xy}

\newcommand{\eg}{{\em e.g.\/}}

\newcommand{\proofbox}{\hfill$\Box$}
\newcommand{\Def}[1]{{\em #1\/}}
\newcommand{\Rem}[1]{{\em #1\/}}
\newcommand{\ex}{\cdot} % explicit multiplication
\newcommand{\ix}{\,} % implicit multiplication
\newcommand{\Nat}{\mathbb{N}}
\newcommand{\equi}{\sim}
\newcommand{\quequi}{=} % equi for use in quotients
\newcommand{\iso}{\cong}

\newcommand{\nonstd}[1]{\underline{#1}}
\newcommand{\Z}{\mathbb{Z}}
\newcommand{\C}{\mathbb{C}}
\newcommand{\rto}{\rightarrow} % reduces to
\newcommand{\rtostar}{\rto^*} % transitive closure of reduces to

\newcommand{\setof}[1]{\{\, #1 \,\}}
\newcommand{\suchthat}{\mid}
\newcommand{\dotminus}{\mbox{$\,-\!\!\!\raisebox{1mm}{$\cdot$}\;$}}
\newcommand{\minusdot}{\dotminus}

\newcommand{\dunion}{\uplus}
\newcommand{\cat}[1]{\mathcal{#1}}
\newcommand{\bref}[1]{(\ref{#1})}

\newtheorem{theorem}{Theorem}
\newtheorem{corollary}[theorem]{Corollary}
\newtheorem{proposition}[theorem]{Proposition}
\newtheorem{example}[theorem]{Example}

\newenvironment{proof}
{\begin{trivlist} \item[] \textsc{Proof:}}{\proofbox \end{trivlist}}
\newenvironment{proof*}
{\begin{trivlist} \item[] \textsc{Proof:}}{\end{trivlist}}
\newenvironment{remark}
{\begin{trivlist}\item[] \textbf{Remark}\: }{\end{trivlist}}

\title
  {An Objective Representation of the Gaussian Integers}

\author
  {Marcelo Fiore\thanks{Research supported by an EPSRC Advanced Research
      Fellowship.}\\
     Computer Laboratory\\ 
     University of Cambridge\\
     United Kingdom%
     \vspace*{2mm}\\
     \texttt{Marcelo.Fiore\mbox{\texttt{@}}cl.cam.ac.uk}
   \and
   Tom Leinster\\
     Department of Mathematics\\ 
     University of Glasgow\\
     United Kingdom%
     \vspace*{2mm}\\
     \texttt{T.Leinster\mbox{\texttt{@}}maths.gla.ac.uk}}

\date{}

\begin{document}

\maketitle

\begin{abstract}
A rig is a ri\Rem{n}g without \Rem{n}egatives.  We analyse the free rig
on a generator $x$ subject to the equivalence $x \equi 1 + x + x^2$,
showing that in it the non-constant polynomials form 
a ring.  This ring can be identified with the Gaussian integers, which
thus acquire objective meaning.
\end{abstract}

\section{Introduction}

Quotient polynomial rings serve as mathematical models in a wide
variety of applications and have been extensively
studied;~see,~\eg,~\cite{GB}.  The corresponding situation for
rigs~(also known as semirings) is underdeveloped. 
The interest for investigating this is that rigs provide direct
mathematical models in scenarios where additive inverses have, 
\Rem{a priori}, no meaning or interpretation.  

One such scenario arises naturally in the context of category
theory, and yields applications in programming and type theory.
Consider the notion of a \Def{distributive category}: a
category with finite sums and finite products with the latter
distributing over the former.  In such a category, sums and products
endow the set of isomorphism classes of objects with the
structure of a rig, its so-called \Def{Burnside rig}.  The Burnside
rig of a distributive category is in fact a ring iff the category is
trivial. 
Thus the natural algebraic structure arising in this context is that
of a rig rather than a ring.  

Following the investigations of Lawvere~\cite{Lawvere} and
Blass~\cite{Blass}, Gates~\cite{Gates} showed that the Burnside rig of
the free distributive category ${\cat D[X]/(X\iso p(X))}$ on a generator
$X$ equipped with an isomorphism $X \iso p(X)$, where $p \in \Nat[x]$
has non-zero constant term, is the quotient polynomial rig
${\Nat[x]/(x\quequi p(x))}$ of the rig 
$\Nat[x]$ under the least congruence identifying $x$ and $p(x)$.  
Thus the structure of $\Nat[x]/(x\quequi p(x))$ and calculations in it
give information on the isomorphisms satisfied by objects 
$X \iso p(X)$ in distributive categories.  For instance, suppose that 
$p_1,p_2\in\Nat[x]$ with $p_1 = p_2$ in $\Nat[x]/(x\quequi p(x))$: then for
all objects $X$ of a distributive category $\cat D$,
\[
X \iso p(X)
\ \implies \  
{p_1(X) \iso p_2(X)}
\quad.
\]
Moreover, every derivation of the equality in the algebra 
$\Nat[x]/(x\quequi p(x))$ yields an isomorphism in the 
category~$\cat D$.    

The distributive categories ${\cat D[X]/(X \iso p(X))}$ 
can be described as categories with objects given by types~(\eg,~the
generator amounts to a recursively defined type) and morphisms given
by programs.  The use of the rig $\Nat[x]/(x\quequi p(x))$ in this
context yields interesting applications to programming and type
theory;~see~\cite{RecTypesIsos} for details. 

\smallskip
In~\cite{FioreLeinster} and~\cite{RecTypesIsos}, we started the study
of the quotient polynomial rigs $\Nat[x]/(x\quequi p(x))$ where 
$p \in \Nat[x]$ has non-zero constant term;~\cite{FioreLeinster}
contains the case of polynomials $p$ with degree at least two,
and~\cite{RecTypesIsos} encompasses all polynomials. 
Among other things, we showed that these quotient polynomial rigs
have a decidable word problem.  The result for polynomials~$p$ of
degree at least two is obtained as a consequence of the following
decomposition:  
\begin{equation}\label{Decomposition}
\Nat[x]/(x \quequi p(x))
\;\; \iso \;\;
\Nat \,\dunion\, \Z[x]/(x-p(x))
\end{equation}
which gives a complete and well-understood description of the rig.  Here
$\dunion$ is disjoint union and the algebraic structure of the right-hand
side 
has additive and multiplicative units respectively given by $0$ and
$1$, addition extended by the obvious action of $\Nat$ on
$\Z[x]/(x-p(x))$, and multiplication extended
freely. 
(The corresponding decomposition result for linear $p$ is more
subtle:~see~\cite{RecTypesIsos}.)

In particular, 
$\Z[x]/(x-p(x))$ embeds as the set of (equivalence classes of)
non-constant polynomials in $\Nat[x]/(x \quequi p(x))$; addition and
multiplication are preserved by this embedding, but the additive and
multiplicative units of $\Z[x]/(x-p(x))$ correspond, inevitably, to
elements of $\Nat[x]/(x \quequi p(x))$ other than $0$ and $1$.  
There are two remarkable aspects to this: first, that the
non-constant elements of the rig $\Nat[x]/(x \quequi p(x))$ carry a
ring structure at all, and second, that this ring is $\Z[x]/(x-p(x))$,
which 
can be thus realised by isomorphism classes of objects in 
${\cat D[X]/(X \iso p(X))}$.   

In this companion paper to~\cite{FioreLeinster,RecTypesIsos} we
analyse one important example of the above situation in detail:
the case $p(x) = 1 + x + x^2$.  There are various reasons
for doing this.  One is that we can establish the
decomposition~\bref{Decomposition} in a very simple, though
insightful, manner, and can prove the further result (akin to
the situation in the theory of Gr\"obner bases for rings) that the
word problem can be solved by a finite strongly normalising
reduction system.  Whether this kind of result holds in generality
is open.  Another motivation, which gives name to the paper, is to
show that the ring of Gaussian integers 
$$\begin{array}{rcl}
\Z[x]/(1+x^2) 
& \iso & 
\Z[i] = \setof{ m + n \ix i \suchthat m,n \in \Z} \subseteq \C
\end{array}$$
has objective meaning in the sense that it arises as the set of
isomorphism classes of objects in a distributive category with the
algebraic operations of addition and multiplication 
corresponding respectively to the categorical operations of sum and
product.  (Recall from above that the additive and multiplicative units cannot
arise as the initial and terminal objects.)  We leave open the 
problem of finding 
a distributive 
category with Burnside rig $\Nat[x]/(x \quequi 1 + x + x^2)$, which would
provide the Gaussian integers with an even more
direct~(\eg,~combinatorial, geometric, or topological) objective
meaning. 

\smallskip
Section~\ref{Results} presents the results of the paper, whilst
Section~\ref{Application} gives an application to programming and 
type theory using the following argument (which we invite the
reader to consider before studying 
the rest of the paper).  Since, as we will see shortly, the identity
${x = x^5}$ holds in ${\Nat[x]/(x\quequi 1 + x + x^2)}$, it follows
that in any distributive category 
the implication   
\begin{equation}\label{X=X5}
X \iso 1 + X + X^2
\ \implies \ 
X \iso X^5
\end{equation}
holds.  In particular, for the distributive category of sets and
functions (with additive structure given by the empty set and
disjoint union, and multiplicative structure given by the singleton
and cartesian product) the set of \Def{Motzkin trees} (that is, unlabelled 
planar unary- and/or binary-branching trees)  clearly satisfies
the hypothesis of the implication~\bref{X=X5}.  Thus, there is an
isomorphism in the language of distributive categories (not merely
in set theory) between the sets of Motzkin trees and five-tuples of
Motzkin trees. 

\section{Results}\label{Results}

A \Def{rig} is a set~$R$ equipped with elements $0$ and $1$ and binary 
operations $+$ and $\ex$ such that $(R,0,+)$ is a commutative monoid, 
$(R,1,\ex)$ is a monoid, and the distributive laws
\begin{center}
\begin{tabular}{ccc}
$0 = a \ix 0$ & \qquad\qquad & $0 = 0 \ix a$ 
\\
$a \ix b + a \ix c = a \ix (b + c)$
& & 
$b \ix a + c \ix a = (b + c)\ix a$
\end{tabular}
\end{center}
hold for all $a, b, c \in R$.

The free rig on a generator $x$ is the set of polynomials $\Nat[x]$
with natural number coefficients equipped with the usual addition and
multiplication of polynomials and their respective units.  The main
object of study in this paper is the quotient polynomial rig
$$
\Nat[x]/(x \quequi 1 + x + x^2)
$$
defined as the quotient rig $\Nat[x]/\!\!\equi$ where $\equi$ is the
smallest congruence on the rig $\Nat[x]$ satisfying 
$x \equi 1 + x + x^2$. 

\medskip
While studying the work of Blass~\cite{Blass} we noticed that there is
an unfolding/folding procedure that works well as a calculational
heuristic method for establishing many identities in quotient
polynomial rigs $\Nat[x]/(x \quequi p(x))$.  In exploring the quotient
polynomial rig ${\Nat[x]/(x \quequi 1 + x + x^2)}$ we soon observed that  
the generator $x$ behaves very much like the imaginary unit.  For
instance, we have 
$x \equi x^5$.  This can be seen from the following calculation,
exemplifying the unfolding/folding procedure referred to above, in
which 
an unfolding step replaces $x^{n+1}$ with 
${x^n + x^{n+1} + x^{n+2}}$~(${n \geq 0}$) whilst a folding step 
does the opposite.
$$\begin{array}[t]{rcll}
x 
& \equi & 1 + x + x^2 
       & \mbox{(unfolding $x$)}
\\[1mm]
& \equi & 1 + x + x + x^2 + x^3  
       & \mbox{(unfolding $x^2$, aiming at cancelling $1$)}
\\[1mm]
& \equi & x + x + x^3  
       & \mbox{(cancelling $1$ and $x^2$ by folding 
                $1 + x + x^2$)}
\\[1mm]
& \equi & x + x + x^2 + x^3 + x^4
       & \mbox{(unfolding $x^3$, aiming at cancelling $x$)} 
\\[1mm]
& \equi & x + x^2 + x^4
       & \mbox{(cancelling $x$ and $x^3$ by folding 
                $x + x^2 + x^3$ )} 
\\[1mm]
& \equi & x + x^2 + x^3 + x^4 + x^5
       & \mbox{(unfolding $x^4$, aiming at cancelling $x$)}
\\[1mm]
& \equi & x^2 + x^4 + x^5
        & \mbox{(cancelling $x$ and $x^3$ by folding 
                $x + x^2 + x^3$)} 
\\[1mm]
& \equi & x^2 + x^3 + x^4 + x^5 + x^5
       & \mbox{(unfolding $x^4$, aiming at cancelling $x^2$)}
\\[1mm]
& \equi & x^3 + x^5 + x^5 
       & \mbox{(cancelling $x^2$ and $x^4$ by folding 
                $x^2 + x^3 + x^4$)}
\\[1mm]
& \equi & x^3 + x^4 + x^5 + x^5 + x^6 
       & \mbox{(unfolding $x^5$, aiming at cancelling $x^3$)}
\\[1mm]
& \equi & x^4 + x^5 + x^6 
       & \mbox{(cancelling $x^3$ and $x^5$ by folding 
                $x^3 + x^4 + x^5$)}
\\[1mm]
& \equi & x^5
       & \mbox{(cancelling $x^4$ and $x^6$ by folding 
                $x^4 + x^5 + x^6$)}
\end{array}$$
The reasons for which this calculation 
goes through 
are explained by the following proposition.

\smallskip
Let $\nonstd{-1} = x^2$, $\nonstd 0 = 1 + \nonstd{-1}$, and 
$\nonstd1 = 1 + \nonstd 0$ in $\Nat[x]$. 
\begin{proposition} \label{Basic_Proposition}
\begin{enumerate}
\item \label{Basic_Proposition_One}
  For $n \geq 0$, $x^n \ix \nonstd 0 \equi \nonstd 0$.

\item \label{Basic_Proposition_Two}
  For all non-constant 
  $p$ in $\Nat[x]$, $p + \nonstd 0 \equi p$.

\item \label{Basic_Proposition_Three}
  For all non-zero 
  $p$ in $\Nat[x]$, $p \ix \nonstd 0 \equi \nonstd 0$.

\item \label{Basic_Proposition_Four}
  For all non-constant 
  $p$ in $\Nat[x]$, $p \ix \nonstd 1 \equi p$.

\item \label{Basic_Proposition_Five}
  For all non-zero 
  $p$ in $\Nat[x]$, $p + \nonstd{-1} \ix p \equi \nonstd 0$.

\item \label{Basic_Proposition_Six}
  For all non-constant 
  $p, q$ in $\Nat[x]$ and for any 
  $r$ in $\Nat[x]$, the \Rem{cancellation law}
  $$
  p + r \equi q + r \implies p \equi q
  $$
  holds.

\item \label{Basic_Proposition_Seven}
  For $p$ in $\Nat[x]$ and $n$ in $\Nat$, $p \equi n$ if and only if 
  $p = n$.
\end{enumerate}
\end{proposition}
\begin{proof}
\bref{Basic_Proposition_One}~$x \ix \nonstd 0 
\ = \ x + x^3
\ \equi \ 1 + x + x^2 + x^3
\ \equi \ 1 + x^2
\ = \ \nonstd 0.$

\bref{Basic_Proposition_Two}~Since $x + \nonstd 0 \equi x$, we also
have that 
$$
x^{n+1} + \nonstd 0 
  \ \equi \ x^{n+1} + x^{n} \ix \nonstd 0
  \    = \ x^n \ix (x + \nonstd 0)
  \ \equi \ x^{n+1}
$$
for all $n \geq 0$.

\bref{Basic_Proposition_Three}~We have
from~\bref{Basic_Proposition_Two} that 
$n \ix \nonstd 0 \ \equi \ \nonstd 0$ for all $n \geq 1$. 
Hence,
$$\begin{array}{rclclcl}
\big(\sum_{i \in I} x^{n_i}\big) \ix \nonstd 0
& = & \sum_{i \in I} (x^{n_i} \ix \nonstd 0)
& \equi & \sum_{i \in I} \nonstd 0
& \equi & \nonstd 0
\end{array}$$
for all finite non-empty~$I$.
  
\bref{Basic_Proposition_Four}~Follows
from~\bref{Basic_Proposition_Two} and~\bref{Basic_Proposition_Three}.

\bref{Basic_Proposition_Five}~Follows
from~\bref{Basic_Proposition_Three}.

\bref{Basic_Proposition_Six}~For $p, q$ non-constant and $r$ non-zero
we have that 
\begin{center}
$p + r \equi q + r 
  \implies p + r + \nonstd{-1} \ix r
  \ \equi \ q + r + \nonstd{-1} \ix r
  \implies p + \nonstd 0  \equi q + \nonstd 0  
  \implies p \equi q
\quad.$
\end{center}

\bref{Basic_Proposition_Seven}~Consider the unique rig homomorphism
from $\Nat[x]/(x \quequi 1 + x + x^2)$ to the rig of countable cardinals
mapping $x$ 
to $\aleph_0$.
\end{proof}
In the light of the proposition, the previous derivation of 
$x \quequi x^5$ in $\Nat[x]/(x=1+x+x^2)$ amounts to the following one:
$$
x 
\ \equi \ x + \nonstd 0 \ix (x^2 + x^3)
\     = \ x + x^2 + x^3 + x^4 + x^5
\     = \ \nonstd 0 \ix (x + x^2) + x^5
\ \equi \ x^5
\quad.
$$

\begin{theorem} \label{Gauss_Theorem}
The subset of 
$\Nat[x]/(x \quequi 1 + x + x^2)$ 
consisting of $($equivalence classes of$)$
non-constant polynomials, equipped
with the usual addition and multiplication but with additive unit 
$\nonstd 0$ 
and multiplicative unit $\nonstd1$, 
is a ring; negatives are given by multiplication with 
$\nonstd{-1}$. 
Further, this ring is $($isomorphic to$)$ the ring of Gaussian integers.
\end{theorem}
\begin{proof}
The first part is a corollary of Proposition~\ref{Basic_Proposition}.  
For the second part, write $R$ for the ring 
in question; then the isomorphism is given by the restriction to $R$ of the 
unique generator-preserving rig homomorphism
$\Nat[x]/(x \quequi 1 + x + x^2) 
   \rightarrow \Z[x]/(1 + x^2)$
and by the unique generator-preserving \emph{ring} homomorphism 
$\Z[x]/(1 + x^2) \rightarrow R$. 
\end{proof}
Explicitly, the isomorphism in the proof amounts to the mappings below.
\begin{center}
$\begin{array}[t]{rcl}
R & \rightarrow & \Z[i] \\
p(x) & \mapsto & p(i) 
\end{array}$
\qquad\qquad\quad
$\begin{array}[t]{rcl}
\Z[i] & \rightarrow & R \\
\begin{array}[t]{c}
\pm m \pm n \ix i \\ 
(m, n \in \Nat) 
\end{array}
& \mapsto & 
\nonstd{\pm1} \ix m + \nonstd{\pm1} \ix n \ix x
\end{array}$
\end{center}
It follows that the Gaussian integers are represented in 
$\Nat[x]/(x \quequi 1 + x + x^2)$ by the 
polynomials 
\begin{equation} \label{GInt_Representation}
\;\;\;\;\;\;\; %%%HACK!!!
m + 1 + x^2 \ \ , \quad
m + n \ix x \ (n \not= 0) \ \ , \quad
m + n \ix x^3 \ (n \not= 0) \ \ , \quad
m \ix x^2 + n \ix x \ \ , \quad
m \ix x^2 + n \ix x^3 
\end{equation}
where $m, n \in \Nat$ are not both $0$. 

\begin{remark}\label{Remark}
Proposition~\ref{Basic_Proposition}\bref{Basic_Proposition_Seven}
and Theorem~\ref{Gauss_Theorem} together imply that 
$\Nat[x]/(x \quequi 1 + x + x^2)$ is formed by extending the addition and
multiplication of the rigs $\Nat$ and $\Z[i]$ to their disjoint union 
$$
\Nat \dunion (\Z \times \Z)
$$
with additive and multiplicative units respectively given by $0$ and
$1$, 
and with addition extended by the obvious action of $\Nat$ on $\Z[i]$:
$$\begin{array}{rclclcl}
\ell + (m,n) 
  & = & (m,n) + \ell
  & = & (\ell + m , n )
    &\quad& (\ell \in \Nat,\, m,n \in \Z) \quad,
\end{array}$$
and multiplication extended freely:
$$\begin{array}{rclclcl}
\ell \ex (m,n) 
  & = & (m,n) \ex \ell 
  & = & \sum_\ell\, (m,n) 
    &\quad& (\ell \in \Nat,\, m,n \in \Z) \quad.
\end{array}$$
\end{remark}

\begin{corollary} \label{Equality_Cor}
For all non-constant 
$p$ and $q$ in $\Nat[x]$ the following are equivalent.
\begin{enumerate}
\item
  $p = q$ in $\Nat[x]/(x \quequi 1 + x + x^2)$.

\item
  $p = q$ in $\Z[x]/(1 + x^2)$.

\item
  $p(i) = q(i)$ in $\Z[i]$.

\end{enumerate}
\end{corollary}
\begin{corollary}
The word problem in $\Nat[x]/(x \quequi 1 + x + x^2)$ is decidable. 
\end{corollary}
\begin{proof}
Given two polynomials in $\Nat[x]$, if they are both non-constant then
evaluate them at $i$ and test for equality in $\Z[i]$; otherwise, by
Proposition~\ref{Basic_Proposition}\bref{Basic_Proposition_Seven}, they
are equivalent if and only if they are equal.
\end{proof}

Our analysis yields an algorithm for obtaining a derivation of the
equality of two polynomials in 
$\Nat[x]/(x \quequi 1 + x + x^2)$.  Indeed, for non-constant 
$p$ and $q$ in $\Nat[x]$ use the division algorithm in $\Z[x]$ to obtain
$$\begin{array}{rcl}
p(x) - q(x) 
& = & 
(w_1(x) - w_2(x)) \ix (1 + x^2) + r(x)
\end{array}$$
with $w_1, w_2$ in $\Nat[x]$ and with remainder $r$ satisfying $r = 0$ or
$0 \leq \deg(r) \leq 1$.  By Corollary~\ref{Equality_Cor}, $p$ and $q$ 
are equal in
$\Nat[x]/(x \quequi 1 + x + x^2)$ if and only if $r = 0$.  In that case we
can obtain a derivation of the equality by noticing that
$$\begin{array}{rcl}
p(x) + (w_1(x) + w_2(x)) \ix x
& \equi &
p(x) + w_1(x) \ix x + w_2(x) \ix (1 + x + x^2) 
\\[1mm]
& = &
q(x) + w_1(x) \ix (1 + x + x^2) + w_2(x) \ix x 
\\[1mm]
& \equi & 
q(x) + (w_1(x) + w_2(x)) \ix x
\end{array}$$
and then deriving $p \equi q$ using the cancellation law
(Proposition~\ref{Basic_Proposition}\bref{Basic_Proposition_Six}). 

\begin{example}\label{2+x2=x4}
Since $2 + i^2 = i^4$ in $\Z[i]$, it follows that 
$2 + x^2 = x^4$ in 
$\Nat[x]/(x \quequi 1 + x + x^2)$.  A 
derivation of this equality using the above method follows. 
$$\begin{array}{rcl}
2 + x^2
& \equi & 2 + x^2
         + (2 + x^2) \ix x
         + \nonstd{-1} \ix (2 + x^2) \ix x
\\[1mm]
& \equi & 2 + x^2
         + 2 \ix x
         + x^2 \ix (1 + x + x^2)
         + \nonstd{-1} \ix (2 + x^2) \ix x
\\[1mm]
&    = & x^4
         + 2 \ix (1 + x + x^2)
         + x^2 \ix x
         + \nonstd{-1} \ix (2 + x^2) \ix x
\\[1mm]
& \equi & x^4
         + (2 + x^2) \ix x
         + \nonstd{-1} \ix (2 + x^2) \ix x
\\[1mm]
& \equi & x^4
\end{array}$$
It is interesting to note that a more direct derivation of the above  
can be obtained by the unfolding/folding procedure:
$$\begin{array}[t]{rcll}
2 + x^2
& \equi & 1 + 1 + x + x^2 + x^3
       & \mbox{$($unfolding $x^2$, aiming at cancelling $1)$}
\\[1mm]
& \equi & 1 + x + x^3
       & \mbox{$($cancelling $1$ and $x^2$ by folding 
                $1 + x + x^2)$}
\\[1mm]
& \equi & 1 + x + x^2 + x^3 + x^4
       & \mbox{$($unfolding $x^3$, aiming at cancelling $1)$}
\\[1mm]
& \equi & x + x^3 + x^4
       & \mbox{$($cancelling $1$ and $x^2$ by folding 
                $1 + x + x^2)$}
\\[1mm]
& \equi & x + x^2 + x^3 + x^4 + x^4
       & \mbox{$($unfolding $x^3$, aiming at cancelling $x)$}
\\[1mm]
& \equi & x^2 + x^4 + x^4
       & \mbox{$($cancelling $x$ and $x^3$ by folding 
                $x + x^2 + x^3)$}
\\[1mm]
& \equi & x^2 + x^3 + x^4 + x^4 + x^5
       & \mbox{$($unfolding $x^4$, aiming at cancelling $x^2)$}
\\[1mm]
& \equi & x^3 + x^4 + x^5
       & \mbox{$($cancelling $x^2$ and $x^4$ by folding 
                $x^2 + x^3 + x^4)$}
\\[1mm]
& \equi & x^4
       & \mbox{$($cancelling $x^3$ and $x^5$ by folding 
                $x^3 + x^4 + x^5)$} 
\end{array}$$
\end{example}

\begin{theorem}
Two polynomials in $\Nat[x]$ are equal in 
$\Nat[x]/(x \quequi 1 + x + x^2)$ if and only if they have the same
normal form in the following strongly normalising reduction system. 
$$\left\{
\begin{array}{rclcl}
x^4 & \rto & 2 + x^2 && 
\\[1mm]
x + x^3 & \rto & 1 + x^2 && 
\\[1mm]
x^n + 1 + x^2 & \rto & x^n && (1 \leq n \leq 3)
\end{array}
\right.$$
\end{theorem}
\begin{proof*}
The reduction system is terminating, as whenever $p \rto q$ we
have that $p(2) > q(2)$.  Further, all critical pairs are joinable 
and so the reduction system is also confluent. 

To conclude the proof we show that the normal forms are exactly 
given by the constants together with the polynomial 
representation~\bref{GInt_Representation} of the Gaussian integers.  
That is, the normal forms are the polynomials 
$$
m + 1 + x^2 \ \ , \quad
m + n \ix x \ \ , \quad
m + n \ix x^3 \ \ , \quad
m \ix x^2 + n \ix x \ \ , \quad
m \ix x^2 + n \ix x^3 
$$
with $m, n \in \Nat$.

By successive applications of the first reduction rule (in the form of 
$x^{m+4} \rto 2 \ix x^m + x^{m+2}$) every polynomial reduces to one of
degree less than or equal to~3.  Further, since  
${a + b \ix x + c \ix x^2 + d \ix x^3}$ ${(a, b, c, d \in \Nat)}$ reduces to
$$
(a+\min(b,d)) + (c+\min(b,d)) \ix x^2 
  + (b \dotminus d) \ix x + (d \dotminus b) \ix x^3
\quad ,
$$
normal forms are either of the form~$(i)~k + \ell \ix x^2 + n \ix x^3$ 
or~$(ii)~k + \ell \ix x^2 + n \ix x$ with $k, \ell, n \in \Nat$.  
We analyse each case in turn. 
\begin{enumerate}
\item[$(i)$] 
  If $n = 0$ then the polynomial is of the form $k + \ell \ix x^2$; in
  which case, if $\ell > k$ it reduces to $(\ell-k) \ix x^2$,
  and if $\ell \leq k$ it reduces to $k$ if $\ell = 0$ and to
  $(k-\ell)+ 1+ x^2$ if $\ell \not= 0$.

  If $n \not= 0$ then the polynomial reduces to 
  $(k \dotminus \ell) +  (\ell \dotminus k) \ix x^2 + n \ix x^3$
  which is either of the form $m + n x^3$ or $m \ix x^2 + n \ix x^3$
  with $m, n \in \Nat$.

\item[$(ii)$] 
  If $n = 0$ then the polynomial is of the form $k + \ell \ix x^2$, and we
  are in the situation of the first case above.

  If $n \not= 0$ then the polynomial reduces to 
  $(k \minusdot \ell) + (\ell \minusdot k) \ix x^2 + n \ix x$ which is
  either of the form $m + n \ix x$ or $m \ix x^2 + n \ix x$ with 
  $m, n \in \Nat$.\proofbox
\end{enumerate}
\end{proof*}
It follows that the word problem in $\Nat[x]/(x \quequi 1 + x + x^2)$ is
decidable in polynomial time.

\begin{example}\label{TwoExamples}
\begin{enumerate}
\item\label{TwoExamplesOne}
For $m \geq 1$, we have that 
$x^{m+4} \rto  2 \ix x^m + x^{m+2} \rto x^m + x^{m-1} + x^{m+1} \rto x^m$.
Hence, as we saw in the introduction, $x^5 = x$ in 
$\Nat[x]/(x \quequi 1 + x + x^2)$. 

\item
In $\Nat[x]/(x \quequi 1 + x + x^2)$, we have that 
$x \ix (1+x^3)^8 \equi 16 \ix x$.  Indeed, 
$(1+x^3)^2 
          = 1 + 2 \ix x^3 + x^6 
   \rtostar 1 + 2 \ix x^3 + x^2 
       \rto 2 \ix x^3$.  
It follows that
$(1+x^3)^4 \equi 4 \ix x^6 \rtostar 4 \ix x^2$ and so
$(1+x^3)^8 \equi 16 \ix x^4 \rtostar 32 + 16 \ix x^2 \rtostar 17 + x^2$.
Finally, 
$x \ix (1+x^3)^8 \equi 17 x + x^3 \rto 16 \ix x + 1 + x^2 \rto 16 \ix x$.
\end{enumerate}
\end{example}

\section{An application}\label{Application}

We conclude the paper with an application to programming and type
theory.

As briefly mentioned in the introduction, the rig 
$\Nat[x]/(x \quequi 1 + x + x^2)$ has straightforward objective
realisation by types; see~\cite{RecTypesIsos} for details.  Indeed,
in the programming language ML, the generator is realised by the
type of Motzkin trees  
defined as follows.
\begin{quote}\begin{verbatim}
datatype  X = e | s of X | m of X * X 
\end{verbatim}\end{quote}
Importantly, calculations in the rig translate as programs that
establish isomorphisms between the associated types.  Thus, for 
instance, the identity $x = x^5$ in 
${\Nat[x]/(x \quequi 1 + x + x^2)}$ entails an isomorphism (in the
language of distributive categories) between Motzkin trees and
five-tuples of Motzkin trees, 
and using the methods of this paper a program realising it can be
automatically constructed.  We illustrate this by working this
example out manually.  

First, consider the second derivation in Example~\ref{2+x2=x4}
establishing the identity $x^4 = 2 + x^2$ in 
${\Nat[x]/(x \quequi 1 + x + x^2)}$.  It yields an isomorphism between
\begin{quote}\begin{verbatim}
type X4 = X * X * X * X 
\end{verbatim}\end{quote}
and 
\begin{quote}\begin{verbatim}
datatype  U = o1 | o2 | p of X * X  
\end{verbatim}\end{quote}
given explicitly by the following program.
\begin{quote}\begin{minipage}{14cm}\begin{verbatim}
val fold1: X4 -> U = fn t => case t of

    ( e, e, e, e )            =>  o1

  | ( e, e, e, s(e) )         =>  o2

  | ( e, e, e, s(s(t)) )      =>  p( e, t )

  | ( e, e, e, s(m(t1,t2)) )  =>  p( s(t1), t2 )

  | ( e, e, e, m(t1,t2) )     =>  p( m(e,t1), t2 )

  | ( e, e, s(t1), t2 )       =>  p( m(s(e),t1), t2 )

  | ( e, e, m(t1,t2), t3 )    =>  p( m(s(s(t1)),t2), t3 )

  | ( e, s(t1), t2, t3 )      =>  p( m(s(m(e,t1)),t2), t3 )

  | ( e, m(t1,t2), t3, t4 )   =>  p( m(s(m(s(t1),t2)),t3), t4 )

  | ( s(t1), t2, t3, t4 )     =>  p( m(m(t1,t2),t3), t4 ) 

  | ( m(t1,t2), t3, t4, t5 )  =>  p( m(s(m(m(t1,t2),t3)),t4), t5 )
\end{verbatim}\end{minipage}\end{quote}
Now, following Example~\ref{TwoExamples}\bref{TwoExamplesOne}, we
exhibit an isomorphism between the types \texttt{X * U} and
\texttt{X}.  A program corresponding to the 
derivation
$$
x \ix (2 + x^2)
\ = \ 2 \ix x + x^3 
\ \equi \ 1 + 2 \ix x + x^2 + x^3
\ \equi \ 1 + x + x^2 
\ \equi \ x
$$
follows.
\begin{quote}\begin{minipage}{14cm}\begin{verbatim}
val fold2: X * U -> X = fn t => case t of

    ( t, o1 )         =>  s(t)

  | ( e, o2 )         =>  e

  | ( s(t), o2 )      =>  m(e,t)

  | ( m(t1,t2), o2 )  =>  m(s(t1),t2)

  | ( t1, p(t2,t3) )  =>  m(m(t1,t2),t3) 
\end{verbatim}\end{minipage}\end{quote}
Finally, an isomorphism between the types \texttt{X * X4} and
\texttt{X} can be given by composing the previous programs: 
\begin{quote}\begin{minipage}{14cm}\begin{verbatim}
val fold: X * X4 -> X = fn t => case t of

    ( t1, t2to5 )  =>  fold2( t1, fold1( t2to5 ) ) 
\end{verbatim}\end{minipage}\end{quote}

\paragraph{Acknowledgements.} 
Our calculations fell into place after a conversation with Bill
\mbox{Lawvere} in which he mentioned a result of Steve Schanuel that
the infinite dimensional elements~(see~\cite{Schanuel}) of some
quotient polynomial rigs actually form a ring.  This led to the
results of this paper and the generalisations presented
elsewhere~\cite{FioreLeinster,RecTypesIsos}.  We also thank David
Corfield for pointing out the relation with Motzkin numbers.

\end{document}